%% file: forArXiv.tex
\documentclass[12pt, fullpage]{article}
\input{./frontstuff}
\usepackage{./mywidefile}

\renewcommand{\baselinestretch}{1.35}
\begin{document}

\title{When is the majority-vote classifier beneficial?}
\author{Mu Zhu\\
University of Waterloo\\
Waterloo, Ontario, Canada N2L 3G1\\
\url{m3zhu@uwaterloo.ca}}
\maketitle

\input{./main}

\section*{Acknowledgments}
The author's research is supported by the Natural Sciences and 
Engineering Research Council (NSERC) of Canada. He would like to thank an Associate Editor and two anonymous referees for their helpful comments.

\bibliographystyle{/u/m3zhu/natbib}
\bibliography{/u/m3zhu/mzstat}

\end{document}

%% file: frontstuff.tex
\usepackage{url} 
\usepackage{float}
\usepackage{graphicx}
\usepackage{amsfonts, amssymb, amscd}
\usepackage{amsmath, delarray, theorem}
\usepackage{natbib}
\bibpunct{(}{)}{;}{a}{}{,}
\renewcommand{\cite}{\citep}

\newcommand{\bitem}{\begin{itemize}}
\newcommand{\eitem}{\end{itemize}}
\newcommand{\benum}{\begin{enumerate}}
\newcommand{\eenum}{\end{enumerate}}
\newcommand{\beqnn}{\begin{eqnarray*}}
\newcommand{\eeqnn}{\end{eqnarray*}}
\newcommand{\beqn}{\begin{eqnarray}}
\newcommand{\eeqn}{\end{eqnarray}}
\newcommand{\normal}[2]{\mbox{N}(#1,#2)}

\newcommand{\E}{\mbox{E}}
\newcommand{\var}{\mbox{Var}}
\newcommand{\cov}{\mbox{Cov}}
\newcommand{\corr}{\mbox{Corr}}

\newcommand{\myvec}[1]{\mathbf{#1}}

\newcommand{\bfourmatrix}{\begin{array}({cccc})}
\newcommand{\efourmatrix}{\end{array}}
\newcommand{\bvector}{\begin{array}({c})}
\newcommand{\evector}{\end{array}}

\newcommand{\btab}{\renewcommand{\baselinestretch}{1}\begin{table}[hptb]}
\newcommand{\etab}{\end{table}\renewcommand{\baselinestretch}{1}}
\newcommand{\bfig}{\renewcommand{\baselinestretch}{1}\begin{figure}[hptb]}
\newcommand{\efig}{\end{figure}\renewcommand{\baselinestretch}{1}}

{\theoremstyle{plain}
  }
{\theoremstyle{plain}
  }
{\theoremstyle{plain}
  }
{\theoremstyle{plain}
  }
{\theoremstyle{plain}
  }
{\theoremstyle{plain}
  \newtheorem{theorem}{Theorem}}
{\theoremstyle{plain}
  }
{\theoremstyle{plain}
  }
{\theoremstyle{plain}
  }
{\theoremstyle{plain}
  }
{\theoremstyle{plain}
  }
{\theorembodyfont{\rmfamily}
  }
{\theorembodyfont{\rmfamily}
  }

%% file: main.tex
\begin{abstract} In his seminal work, \citet{weak-learner} proved 
that weak classifiers could be improved to achieve arbitrarily high 
accuracy, but he never implied that a simple majority-vote mechanism could 
always do the trick. By comparing the asymptotic misclassification error of the majority-vote classifier with the average individual error, we discover an interesting phase-transition phenomenon. For binary classification with equal prior probabilities, our result implies that, for the majority-vote mechanism to work, the 
collection of weak classifiers must meet the minimum requirement of having an average true 
positive rate of at least $50\%$ and an average false positive rate of at most 
$50\%$. \end{abstract}

\vspace{5mm}
{\bf Key words}: bagging; boosting; ensemble learning; phase transition; 
random forest; weak learner.

\pagebreak

\def\pr{\mbox{Pr}}

\section{Introduction}

Consider the binary classification problem. Let $y\in\{0,1\}$ be the class 
label, with $\pi=\pr(y=1)$ being the prior probability of class 1. Let 
$\mathcal{F}_n = \{f_1, f_2, ..., f_n\}$ be a collection of binary classifiers. 
We are interested in the so-called majority-vote classifier, 
\[
 g_n = \sum_{i=1}^n f_i.
\]
For simplicity, 
we assume that a strict majority is needed for class 1, i.e., 
\[
 \hat{y} = \begin{cases}
 1, & \mbox{if} \quad g_n > n/2; \\
 0, & \mbox{if} \quad g_n \leq n/2.
 \end{cases}
\]
Classifiers that perform slightly better than random guessing are referred 
to as ``weak learners'' \citep{weak-learner-orig, weak-learner}. It is 
widely held \citep[e.g.,][]{webb-pattern-book, kuncheva-book,
ieee-majvote} that, as long as a collection of weak learners 
$f_1, f_2, ..., f_n$ are {\em nearly} independent of each other, the majority-vote classifier $g_n$ will achieve greater accuracy.

The purpose of this article is to show that, in order for $g_n$ to have better 
performance than a typical individual $f_i$, the collection of weak learners $\mathcal{F}_n$ cannot be ``too weak'', \emph{even if the total amount of correlation among them is well under control}. For example, assuming 
equal priors ($\pi=1/2$), the collection must have an average true positive rate 
of at least $50\%$ and an average false positive rate of at most $50\%$.

\subsection{Common misconceptions}
\label{sec:wrong}

In his seminal paper, \citet{weak-learner} proved that weak 
learning algorithms --- ones that perform only slightly better than random 
guessing --- can be turned into ones capable of achieving arbitrarily high 
accuracy. Does our conclusion stand in the face of Schapire's result? 
Certainly not. Schapire's constructive proof shows how one can {\em 
sequentially} boost weak learners into strong ones 
\citep[see also][]{boosting-orig}, but his proof 
certainly does {\em not} imply that one can always do so with a simple 
majority-vote mechanism --- a common misconception. 

Such widespread misconceptions are perhaps fueled partially by incomplete 
understandings of Breiman's influential algorithms --- namely, bagging 
\citep{bagging} and random forest \citep{randomForest} --- that do indeed 
use the majority-vote mechanism, and partially by the popular lessons 
drawn from the highly-publicized million-dollar Netflix contest 
\cite[see, e.g.,][]{netflix-wired-news} that testify to the wisdom of 
crowds \citep{wisdom-crowds}. But it is a logic flaw to conclude that the 
majority-vote mechanism can be used to improve {\em any} collection of 
weak learners simply because it has had a number of successful 
applications. Even for \citet{wisdom-crowds}, crowd wisdom is not always 
the result of simple majority voting. An important goal that we hope to 
accomplish with our analysis is to put a halt to the spread of such common 
misconceptions.

\section{Analysis}
\label{sec:main}
\def\err{\mbox{err}}
\def\Err{\mbox{Err}}

\begin{figure}[htpb]
\centering
\includegraphics[width=0.618\textwidth, angle=270]{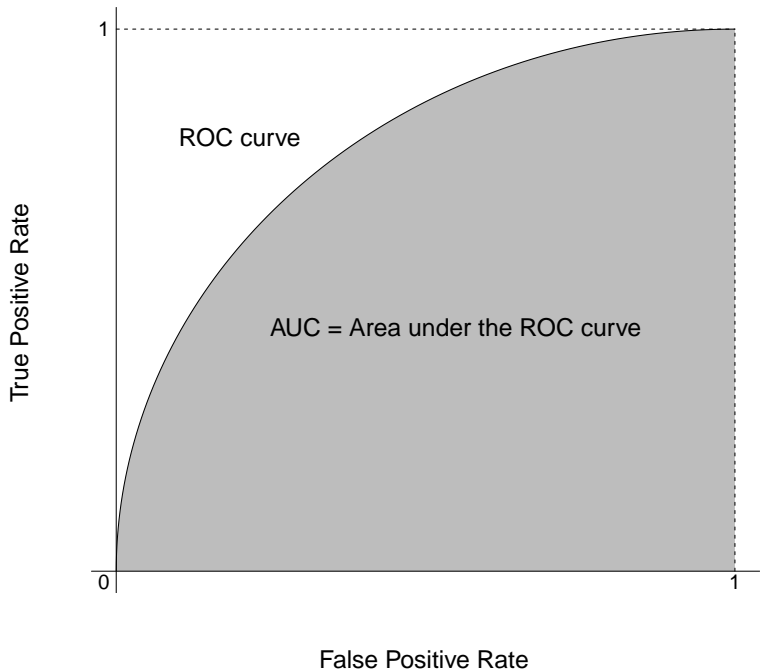}
\caption{\label{fig:ROC}%
An illustration of the ROC curve.}
\end{figure}

To make our point, let
\beqn
\label{eq:TP}
 p_i \equiv \pr(f_i=1|y=1) 
\eeqn
denote the true positive rate (TPR) of the classifier $f_i$ and
\beqn
\label{eq:FP}
 q_i \equiv \pr(f_i=1|y=0), 
\eeqn
its false positive rate (FPR). 
These two quantities form the key signatures of a binary classifier. For example, the widely-used receiver-operating characteristic (ROC) curve \citep[e.g.,][]{roc-early1, roc-early2, roc-popularizer, roc-science, rocbk} has the FPR on the horizontal axis and the TPR on the vertical axis (Figure~\ref{fig:ROC}).

Algorithms that use the majority-vote mechanism such as the random forest \citep{randomForest} often employ an i.i.d.~mechanism to generate 
the collection 
\beqn
\label{eq:RF}
 \mathcal{F}_n = \{
 f_i(\myvec{x}) = f(\myvec{x}; \theta_i)
 \quad\mbox{where}\quad
 \theta_i \overset{iid}{\sim} \mathcal{P}_{\theta}\}.
\eeqn
Here, $f(\myvec{x};\theta_i)$ is a classifier completely parameterized by $\theta_i$, and the statement ``$\theta_i \overset{iid}{\sim} \mathcal{P}_{\theta}$'' means that each $f(\myvec{x}; \theta_i)$ is generated using an i.i.d. stochastic mechanism, $\mathcal{P}_{\theta}$, for example, bootstrap sampling \citep[e.g.,][]{bagging} and/or random subspaces \citep[e.g.,][]{tkho}.

Therefore, we can define
\beqn
\label{eq:p-def}
 p \equiv \E_{\mathcal{P}_{\theta}}(p_i) = \int_{\theta_i} \pr(f(\myvec{x}; \theta_i)=1|y=1) d\mathcal{P}_{\theta} 
\eeqn
and
\beqn
\label{eq:q-def}
 q \equiv \E_{\mathcal{P}_{\theta}}(q_i) = \int_{\theta_i} \pr(f(\myvec{x}; \theta_i)=1|y=0) d\mathcal{P}_{\theta} 
\eeqn
to be the {\em average} TPR and the {\em average} FPR of the collection $\mathcal{F}_n$, respectively.
We shall further assume that both $p$ and $q$ are on the {\em open} interval of $(0,1)$, ignoring the trivial cases of $p, q \in \{0,1\}$. 

\subsection{Strategy}
\label{sec:strategy}

We will write
\beqnn
 \mu_p \equiv \E_p(g_n) \equiv \E(g_n|y=1), &\quad& \var_p(g_n) \equiv \var(g_n|y=1), \\
 \mu_q \equiv \E_q(g_n) \equiv \E(g_n|y=0), &\quad& \var_q(g_n) \equiv \var(g_n|y=0).
\eeqnn
Let
\beqn
\label{eq:asymp-var}
\sigma_p^2 \equiv \lim_{n\rightarrow\infty}\left[\frac{\var_p(g_n)}{n}\right] \quad\mbox{and}\quad 
\sigma_q^2 \equiv \lim_{n\rightarrow\infty}\left[\frac{\var_q(g_n)}{n}\right].
\eeqn
Our argument is asymptotic for $n\rightarrow\infty$. That is to say, we 
assume that we have access to a fairly large collection of classifiers. 
This is certainly the case for both bagging \citep{bagging} and random 
forest \citep{randomForest}. Our strategy is very simple. The error rate of the individual classifier $f_i$ is given by
\beqnn
 \err_i 
&\equiv& \pr(f(\myvec{x};\theta_i)=0|y=1) \pr(y=1) + \pr(f(\myvec{x};\theta_i)=1|y=0) \pr(y=0) \\
&=& (1 - p_i) \pi + q_i (1 - \pi).
\eeqnn
Thus, the {\em mean} individual error of classifiers in the collection $\mathcal{F}_n$ is
\beqn
\err 
\equiv \E_{\mathcal{P}_{\theta}} ( \err_i ) 
= \E_{\mathcal{P}_{\theta}} \left[ (1-p_i) \pi + q_i (1-\pi) \right] 
= (1-p) \pi + q (1-\pi), \label{eq:err} 
\eeqn 
whereas the error rate of the majority-vote classifier, $g_n$, is simply
\beqnn
\Err(n)&\equiv&
 \pr\left(g_n \leq n/2 \big | y=1\right) \pr(y=1) + 
 \pr\left(g_n > n/2 \big | y=0\right) \pr(y=0).
\eeqnn
First, we apply the central limit theorem (some technical details in Section~\ref{sec:tech})
and estimate $\Err(n)$ by
\beqn
\label{eq:hatErr}
\widehat{\Err}(n) = 
  \left[ 
    \Phi\left( \frac{n/2-\mu_p}{\sigma_p\sqrt{n}} \right)
  \right] \pi +
  \left[
  1-\Phi\left( \frac{n/2-\mu_q}{\sigma_q\sqrt{n}} \right)
  \right]  (1-\pi),
\eeqn
where $\Phi(\cdot)$ is the cumulative distribution function of the 
standard $\normal{0}{1}$ random variable.
Then, we compare $\widehat{\Err}(n)$ with the mean individual error, given in (\ref{eq:err}).
The difference, 
\[
 \Delta(n) \equiv \widehat{\Err}(n)-\err,
\] 
is one way to measure whether $g_n$ provides any genuine 
benefit over a typical $f_i$, with $\Delta(n) < 0$ being an indication of true 
benefit. We pay particular attention to the limiting quantity,
\[
 \Delta(\infty) \equiv \lim_{n\rightarrow\infty} \Delta(n),
\]
which turns out to have interesting phase transitions in the space of 
$(p,q)$; see Section~\ref{sec:phase} below.

\subsection{Details}
\label{sec:tech}

In order to apply the central limit theorem so as to estimate $\Err(n)$ by (\ref{eq:hatErr}), the key technical requirement is that both limits $\sigma_p$ and $\sigma_q$ are finite and non-zero \citep[see, e.g.,][]{billingsley, mixingCLT-ref}.


\begin{theorem}
\label{thm:main}
If
\beqn
\label{eq:proper-condition}
0 < \sigma_p, \sigma_q < \infty,
\eeqn
then
\[
 \lim_{n\rightarrow\infty} \Err(n) =
 \lim_{n\rightarrow\infty} \widehat{\Err}(n),
\]
and the limit (given in
Table~\ref{tab:limErr}) has phase transitions in the space of $(p,q)$
depending on whether $p$ and $q$ are larger 
than, equal to, or less than $1/2$.
\end{theorem}
Using the result of Theorem~\ref{thm:main} and
Table~\ref{tab:limErr} in particular, 
we can compute easily --- for any given position $(p,q)$ in the open
unit square --- the value of
\[
 \Delta(\infty) 
 = \lim_{n\rightarrow\infty}\widehat{\Err}(n) - \err,
\] 
and determine whether $\Delta(\infty)$ is positive, negative, or zero. The 
proof of Theorem~\ref{thm:main} given below (Section~\ref{sec:proof}) 
explains
the origin of the phase-transition phenomenon, but 
readers not interested in the details may skip it and jump right to 
Section~\ref{sec:phase}. 

\begin{table}[htpb]
\caption{\label{tab:limErr}%
The limiting behavior of $\widehat{\Err}(n)$ as $n \rightarrow \infty$.} 
\centering
\fbox{
\begin{tabular}{l|ccc}
        & $p<1/2$ & $p=1/2$ & $p>1/2$ \\
\hline
$q>1/2$ & $1$         & $1-\pi/2$ & $1-\pi$     \\
$q=1/2$ & $(1+\pi)/2$ & $1/2$     & $(1-\pi)/2$ \\
$q<1/2$ & $\pi$       & $\pi/2$   & $0$         
\end{tabular}}
\end{table}

\subsubsection{Proof of Theorem~\ref{thm:main}}
\label{sec:proof}

The condition $0 < \sigma_p, \sigma_q < \infty$ allows us to apply the central limit theorem for dependent random variables \citep[see, e.g.,][]{billingsley, mixingCLT-ref}, so that
\[
 \frac{g_n-\mu_p}{\sigma_p\sqrt{n}} \bigg | y=1 \overset{D}{\longrightarrow} \normal{0}{1}\quad\mbox{and}\quad
 \frac{g_n-\mu_q}{\sigma_q\sqrt{n}} \bigg | y=0 \overset{D}{\longrightarrow} \normal{0}{1}.
\]
Therefore, 
\[
 \lim_{n\rightarrow\infty}\Err(n) 
=\lim_{n\rightarrow\infty}\widehat{\Err}(n), 
\]
where $\widehat{\Err}(n)$ is given by (\ref{eq:hatErr}). 
Under (\ref{eq:TP})-(\ref{eq:FP}) and (\ref{eq:p-def})-(\ref{eq:q-def}), we have
\beqnn
 \mu_p = \E_p(g_n) 
= \sum_{i=1}^n \E(f_i|y=1)
= \sum_{i=1}^n \int_{\theta_i} \pr(f(\myvec{x};\theta_i)=1|y=1) d\mathcal{P}_{\theta} 
= np.
\eeqnn
Similarly, $\mu_q=\E_q(g_n)=nq$. So, 
\beqn
\label{eq:hatErr-propercase}
\widehat{\Err}(n) 
&=&
  \left[ 
    \Phi\left( \frac{\sqrt{n}(1/2-p)}{\sigma_p} \right)
  \right]  \pi    +
  \left[
  1-\Phi\left( \frac{\sqrt{n}(1/2-q)}{\sigma_q} \right)
  \right]  (1-\pi).
\eeqn
The theorem is proved 
by applying the following limit 
to equation~(\ref{eq:hatErr-propercase}): 
\beqnn
 \lim_{n\rightarrow\infty} 
 \Phi\left( c \sqrt{n} \right) 
&=&
\begin{cases}
1,   & c>0; \\
1/2, & c=0; \\
0,   & c<0.
\end{cases}
\eeqnn


\subsection{Phase transition}
\label{sec:phase}

Table~\ref{tab:limDelta} shows the values of $\Delta(\infty)$ as they vary in the space of $(p,q)$. Notice the abrupt phase transitions occurring at the boundaries, $p=1/2$ and $q=1/2$. In particular, for fixed $q$, there is a jump of size $\pi/2$ going from $p<1/2$, $p=1/2$, to $p>1/2$; for fixed $p$, there is a jump of size $1/2-\pi/2$ going from $q<1/2$, $q=1/2$, to $q>1/2$. 

\def\A{A^{\pi}_{p,q}}
\begin{table}[htpb]
\caption{\label{tab:limDelta}%
The limiting behavior of $\Delta(n)$ as $n \rightarrow \infty$.}
\centering
\fbox{
\begin{tabular}{l|lll}
        & \multicolumn{1}{c}{$p<1/2$} 
        & \multicolumn{1}{c}{$p=1/2$} 
        & \multicolumn{1}{c}{$p>1/2$} \\
\hline
$q>1/2$ & $\A+(1-\pi)$                     & $\A+(1-\frac{3\pi}{2})$ & $\A+(1-2\pi)$     \\
$q=1/2$ & $\A+(\frac{1}{2}-\frac{\pi}{2})$ & $\A+(\frac{1}{2}-\pi)$  & $\A+(\frac{1}{2}-\frac{3\pi}{2})$ \\
$q<1/2$ & $\A$                             & $\A-\frac{\pi}{2}$      & $\A-\pi$          \\
\end{tabular}}\\
{\footnotesize Notation: $\A \equiv \pi p-(1-\pi)q$.}
\end{table}

\begin{figure}[htpb]
\centering
\includegraphics[width=0.325\textwidth, angle=270]{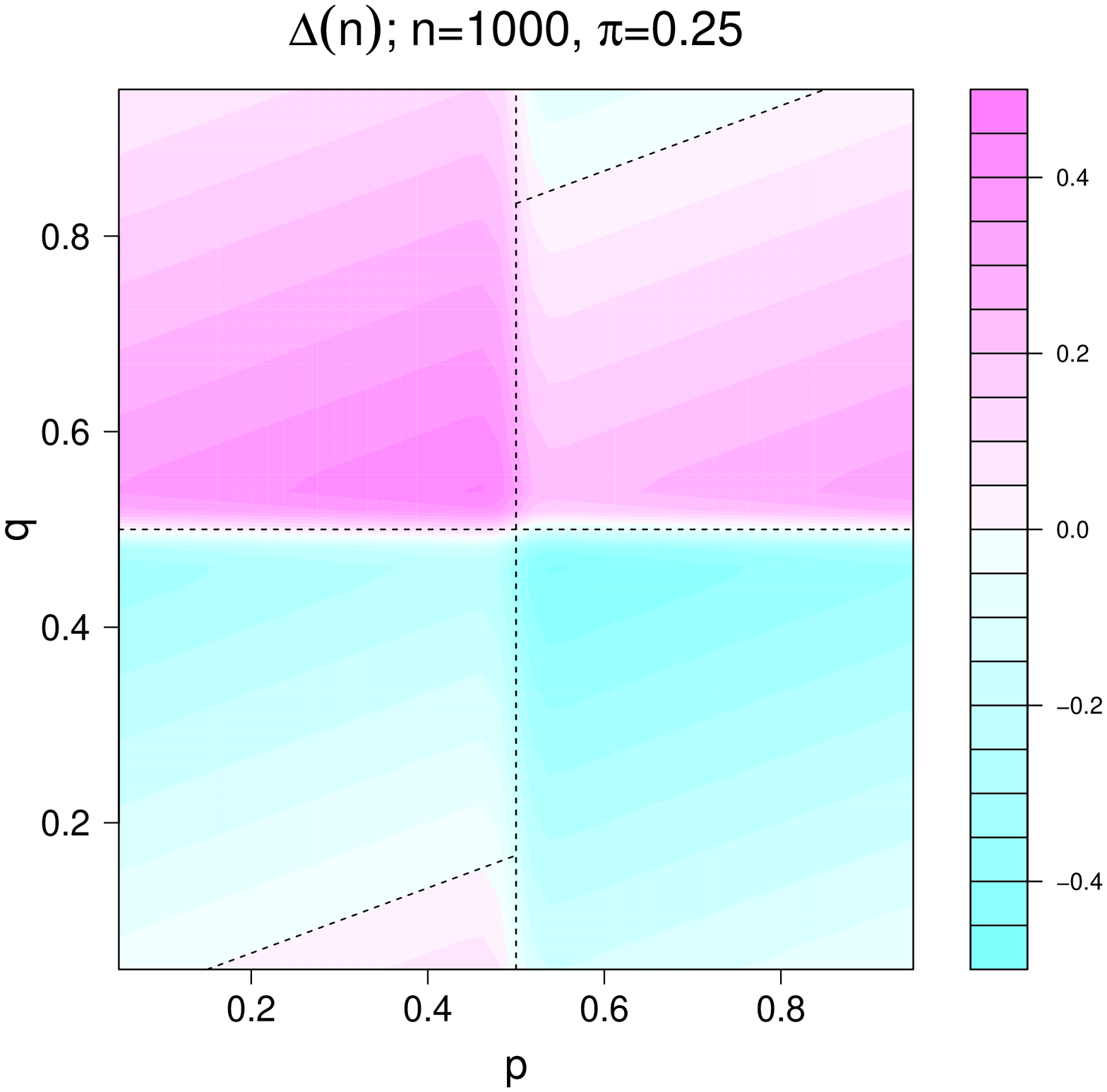}
\includegraphics[width=0.325\textwidth, angle=270]{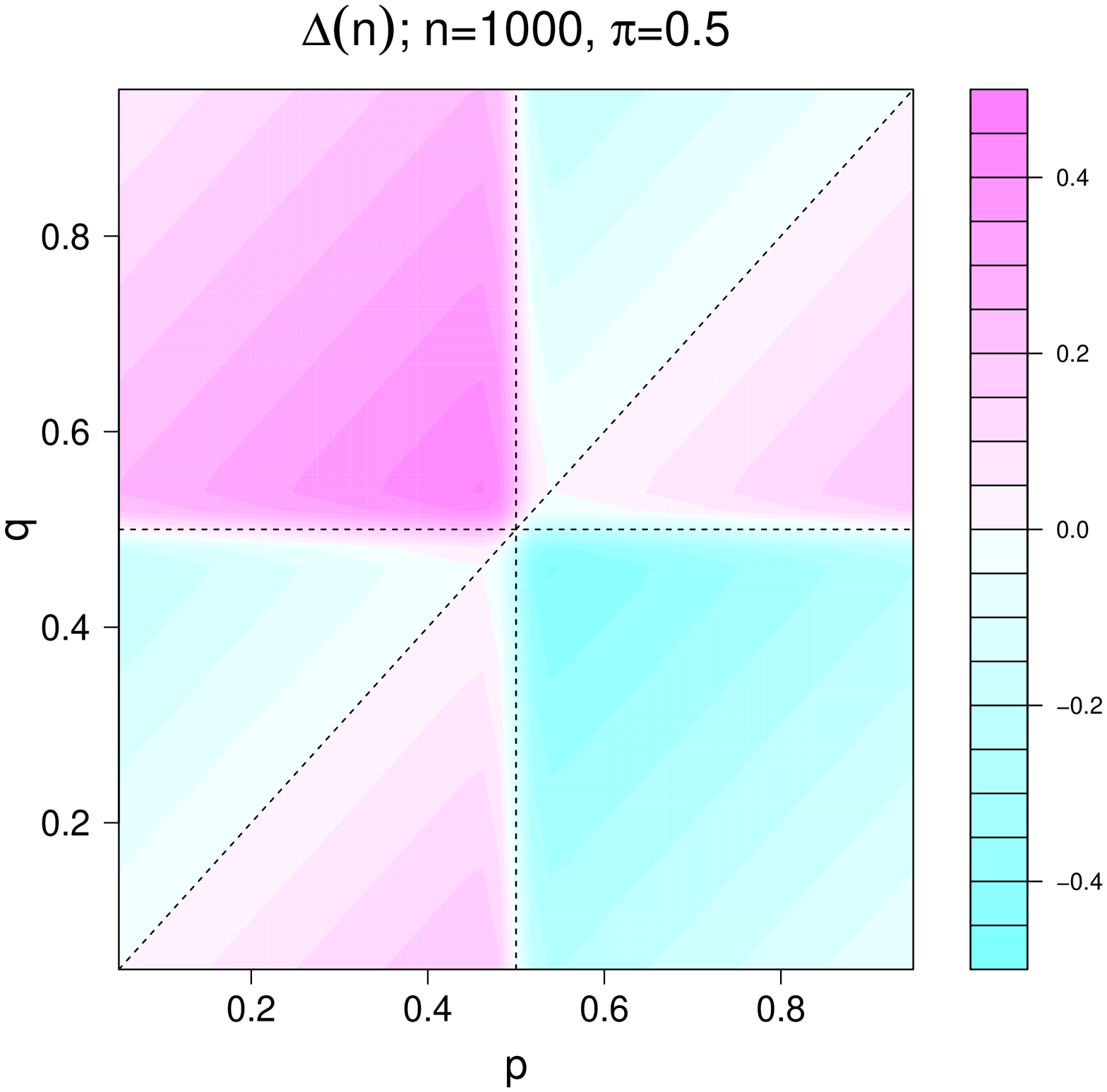} 
\includegraphics[width=0.325\textwidth, angle=270]{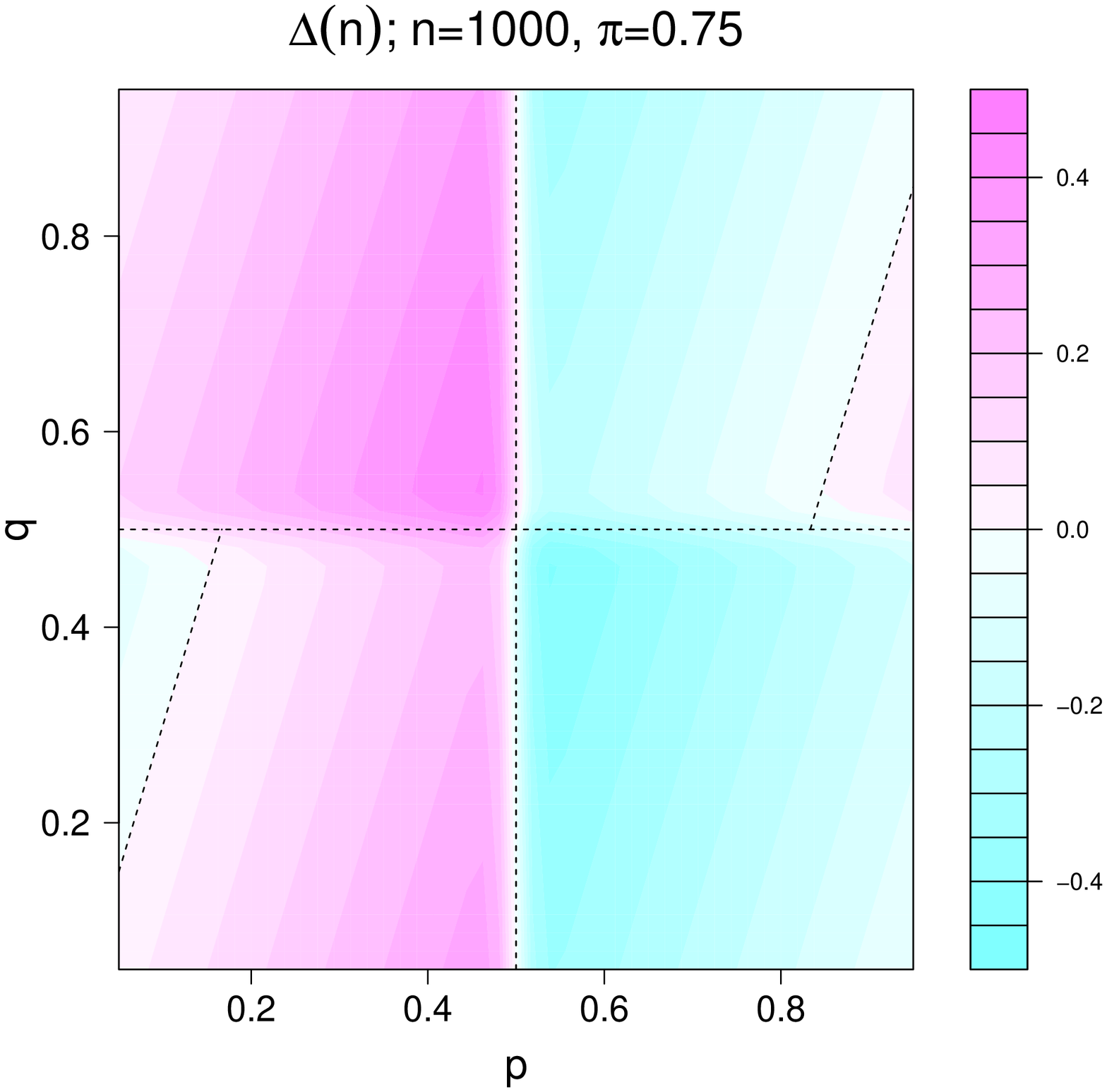}
\caption{\label{fig:pi}%
$\Delta(n)$ as a function of $(p,q)$ for $\pi=0.25$ (left), $\pi=0.5$ (middle), and $\pi=0.75$ (right), while fixing $n=1000$.}
\end{figure}

Figure~\ref{fig:pi} shows $\Delta(n)$ as a function of $(p,q)$ for $\pi=0.25, 0.5$ and $0.75$ while fixing $n=1000$. It is evident from both Figure~\ref{fig:pi} and Table~\ref{tab:limDelta} that $\pi \in (0,1)$ merely changes the shape of the function $\Delta(\infty)$ within each region --- as specified by whether each of $p, q$ is larger than, equal to, or smaller than $1/2$ --- and the size of the jumps between these regions, but it does not cause these phase transitions to appear or to disappear. In other words, the phase transition phenomenon at $p=1/2$ and at $q=1/2$ is universal regardless of what value $\pi$ takes on. 

While a phase diagram can be produced easily for any value of $\pi$, for the remainder of the paper we shall focus on $\pi=1/2$, a case that deserves special attention for binary classification (Table~\ref{tab:limDelta-half}).

\begin{table}[htpb]
\caption{\label{tab:limDelta-half}%
The limiting behavior of $\Delta(n)$ as $n \rightarrow \infty$, for the special case of $\pi=1/2$.}
\centering
\fbox{
\begin{tabular}{l|lll}
        & \multicolumn{1}{c}{$p<1/2$} 
        & \multicolumn{1}{c}{$p=1/2$} 
        & \multicolumn{1}{c}{$p>1/2$} \\
\hline
$q>1/2$ & $(p-q)/2+1/2$ & $(p-q)/2+1/4$ & $(p-q)/2$     \\
$q=1/2$ & $(p-q)/2+1/4$ & $(p-q)/2$     & $(p-q)/2-1/4$ \\
$q<1/2$ & $(p-q)/2$     & $(p-q)/2-1/4$ & $(p-q)/2-1/2$ \\ 
\end{tabular}}
\end{table}

\subsection{The case of $\pi=1/2$}

Recall from (\ref{eq:p-def})-(\ref{eq:q-def}) that $p$ and $q$ are the average TPR and the average FPR of the collection $\mathcal{F}_n$. Intuitively, it is natural for us to require that, on average, individual classifiers in $\mathcal{F}_n$ should be more likely to predict $\hat{y}=1$ when $y=1$ than when $y=0$, that is, we would like to have
\[
p>q.
\] 
Indeed, when $\pi=1/2$, we easily can see that the mean individual error is smaller than the error of random guessing if and only if $p>q$:
\[
\err=(1-p)\pi+q(1-\pi)=\frac{1-(p-q)}{2} < \frac{1}{2} \quad \Leftrightarrow \quad p>q.
\]
Therefore, based on ``conventional wisdom'' (Section~\ref{sec:wrong}), one may be tempted to think that, on average, $\mathcal{F}_n$ can be considered a collection of weak classifiers so long as $p>q$, and that taking majority votes over such a collection should be beneficial. But, according to our preceding analysis, the condition $p>q$ alone is actually ``too weak'' and {\em not} enough by itself to ensure that $g_n$ is truly beneficial.

\begin{figure}[htpb]
\centering
\includegraphics[width=0.5\textwidth, angle=270]{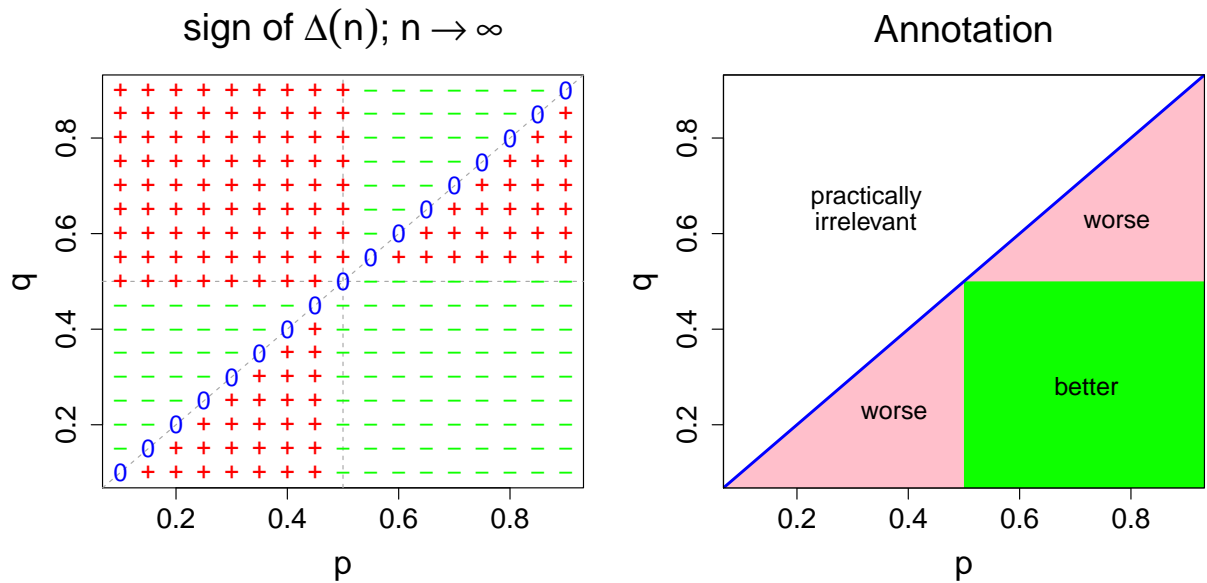}
\caption{\label{fig:phase}%
The phase diagram showing whether $\Delta(\infty)$ is positive ($+$, red), 
negative ($-$, green), 
or zero ($0$, blue), in the space of $(p,q)$, assuming that $\pi=1/2$. 
Notice that abrupt phase transitions 
occur at $p=1/2$ and $q=1/2$.}
\end{figure}

Figure~\ref{fig:phase} contains the corresponding phase diagram showing the sign of $\Delta(\infty)$ for the case of $\pi=1/2$. The phase diagram clearly shows that we need more than $p>q$ in order to have $\Delta(\infty) \leq 0$. More specifically, we see that we need $p$ and $q$ to be on {\em different sides} of $1/2$, that is,
\[
 p \geq 1/2 \geq q.
\]
The condition $p>q$ by itself is too weak. If, for example, 
\[
 p > q > 1/2 \quad\mbox{or}\quad 1/2 > p > q,
\]
then the majority-vote classifier $g_n$ actually performs worse! In other words, even if $p>q$ so that a typical classifier in $\mathcal{F}_n$ can be considered a weak learner, the 
majority-vote mechanism only makes matters worse if $p$ and $q$ are on the {\em same side} of $1/2$. 

This is because, if $p$ and $q$ are on the same side of $1/2$, then, as we take majority votes over more and more weak learners, $g_n$ will eventually classify everything into the same class with probability one (and become useless as a result). In particular, $g_n$ will classify everything into class 1 if both $p$ and $q$ are greater than $1/2$, and into class 0 if both are less than $1/2$.

Perhaps somewhat curiously, symmetry leads to the conclusion that it is possible to have $\Delta(\infty)<0$ even when $p<q$, e.g., if $1/2 < p < q$ or $p < q < 1/2$ --- again, see 
Figure~\ref{fig:phase}. However, this is not so much because $g_n$ is a good classifier, but because the collection $\mathcal{F}_n$ is {\em very} poor when $p<q$; in fact, a typical classifier from $\mathcal{F}_n$ is worse than random, and more likely to predict $\hat{y}=1$ when $y=0$ and vice versa. But, as we've pointed out earlier, if $p$ and $q$ are on the same side of $1/2$, $g_n$ will eventually classify everything into the same class. Now, that's useless indeed, but it surely is still better than classifying things into the opposite class! Though conceptually curious, these cases are not of practical interest. 

\subsection{Conclusion}
\label{sec:interpretation}

What does all this mean? Given a collection of classifiers, $\mathcal{F}_n$, intuitively one may be tempted to think of it as a collection of weak learners if its average TPR is greater than its average FPR (i.e., if $p>q$), but taking majority votes over such a collection cannot be guaranteed to produce improved results. To obtain improvements (assuming equal priors), the collection must have an average TPR of at least $50\%$ (i.e., $p \geq 0.5$) {\em and} an average FPR of at most $50\%$ (i.e., $q \leq 0.5$). In other words, the typical classifier in the collection cannot be arbitrarily weak.

\section{Discussions}
\label{sec:examples}

How realistic is the technical condition of Theorem~\ref{thm:main}, and 
hence how relevant is our conclusion above? In this section, we describe 
two scenarios that satisfy the theorem's condition (Sections~\ref{sec:indep} and \ref{sec:AR}) and one that does not (Section~\ref{sec:equi}).

\subsection{$\var(f_i|\theta_i; y=1)$ versus $\var(f_i|y=1)$}

By the law of total variance, we have
\beqnn
 \var_p(f_i) 
&\equiv& \var(f_i|y=1) \\
&=& 
 \E_{\mathcal{P}_{\theta}}   \left[\var(f_i|\theta_i; y=1)\right] +
 \var_{\mathcal{P}_{\theta}} \left[\E(f_i|\theta_i;y=1)\right] \\
&=&
 \E_{\mathcal{P}_{\theta}}   \left[p_i(1-p_i)\right] +
 \var_{\mathcal{P}_{\theta}} \left[p_i\right] \\
&=&
 \E_{\mathcal{P}_{\theta}}   (p_i) - \E_{\mathcal{P}_{\theta}} (p_i^2) +
 \var_{\mathcal{P}_{\theta}} \left(p_i\right) \\
&=&
 \E_{\mathcal{P}_{\theta}}   (p_i) - \left[ \E^2_{\mathcal{P}_{\theta}} (p_i) + \var_{\mathcal{P}_{\theta}} (p_i) \right] +
 \var_{\mathcal{P}_{\theta}} \left(p_i\right) \\
 &=&
 p(1-p).
\eeqnn
While it is obvious that $\var(f_i|\theta_i;y=1) = p_i(1-p_i)$, it is perhaps not so obvious that, under $\theta_i \overset{iid}{\sim} \mathcal{P}_{\theta}$, 
we simply have $\var(f_i|y=1)=p(1-p)$ as well. The same argument establishes that $\var_q(f_i) = q(1-q)$. Furthermore, following the definitions of $\var_p(f_i)$ and $\var_q(f_i)$, we define 
\[
 \cov_p(f_i, f_j) = \cov(f_i, f_j | y=1), \qquad
 \cov_q(f_i, f_j) = \cov(f_i, f_j | y=0),
\]
as well as 
\[
 \corr_u(f_i,f_j) = \frac{\cov_u(f_i, f_j)}{\sqrt{\var_u(f_i)\var_u(f_j)}}
\]
for $u=p,q$.


\subsection{Case 1} 
\label{sec:indep}

Suppose that $\corr_u(f_i, f_j)=0$ for all $i \neq j$ and $u=p,q$ --- e.g., if $f_1, f_2, ..., f_n$ are conditionally independent given either $y=0$ or $y=1$, then 
\beqnn
\sigma_p^2
=\lim_{n\rightarrow\infty}\left[\frac{\var_p(g_n)}{n}\right]
=\lim_{n\rightarrow\infty}\left[\frac{np(1-p)}{n}\right]
=p(1-p) \in (0, \infty).
\eeqnn
Similarly, we can conclude that $0 < \sigma_q^2 < \infty$ as well, so the condition of Theorem~\ref{thm:main} is satisfied.
In reality, however, the conditional independence assumption is not very 
realistic. 
Since the individual classifiers $f_1, f_2, ..., f_n$ typically 
rely on the same training data $\{(y_i, \myvec{x}_i)\}_{i=1}^N$ to ``go 
after'' the same quantity of interest, namely $\pr(y|\myvec{x})$, there is 
usually an non-ignorable amount of correlation among them.

\subsection{Case 2} 
\label{sec:AR}
\def\ARfact{\frac{1+\gamma}{1-\gamma}}

Alternatively, suppose that, for all $i, j$ and $u=p,q$, $\corr_u(f_i, f_j) = \gamma^{|i-j|}$ for some $0< \gamma <1$. The correct interpretation of this assumption is as follows: between a given classifier $f_i$ and other members of $\mathcal{F}_n$, the largest correlation is $\gamma$; the second largest correlation is $\gamma^2$; and so on. The exponential decay structure is admittedly artificial, but such a structure will ensure that both limits $\sigma_p$ and $\sigma_q$ are, again, finite and non-zero. In particular, 
\beqn
\var_p(g_n) 
&=& \sum_{i=1}^n \var_p(f_i) + \sum_{i \neq j} \cov_p(f_i,f_j) \notag \\
&=& np(1-p) +  
 2\left[(n-1)\gamma + (n-2)\gamma^2 + ... +\gamma^{n-1}\right]p(1-p) \notag \\
&=& np(1-p) \left[1 + 
 2\sum_{j=1}^{n-1} \left(1-\frac{j}{n}\right)\gamma^j 
 \right] \label{eq:liu-preapprox},
\eeqn
which means
\beqn
\sigma_p^2 =
\lim_{n\rightarrow\infty}\left[\frac{\var_p(g_n)}{n}\right]
= p(1-p) \left[1 + 2 \sum_{j=1}^{\infty} \gamma^j \right] 
= p(1-p)\left[\frac{1+\gamma}{1-\gamma}\right] \in (0,\infty) 
 \label{eq:liu-postapprox} 
\eeqn
and likewise for $\sigma_q^2$. 
The limit going from 
(\ref{eq:liu-preapprox}) to (\ref{eq:liu-postapprox}) as $n\rightarrow\infty$ is 
standard and frequently mentioned in the Markov chain Monte Carlo (MCMC) 
literature \citep[e.g.,][Section 5.8]{liu}.

\begin{figure}[htpb]
\centering
\includegraphics[width=0.4\textwidth, angle=270]{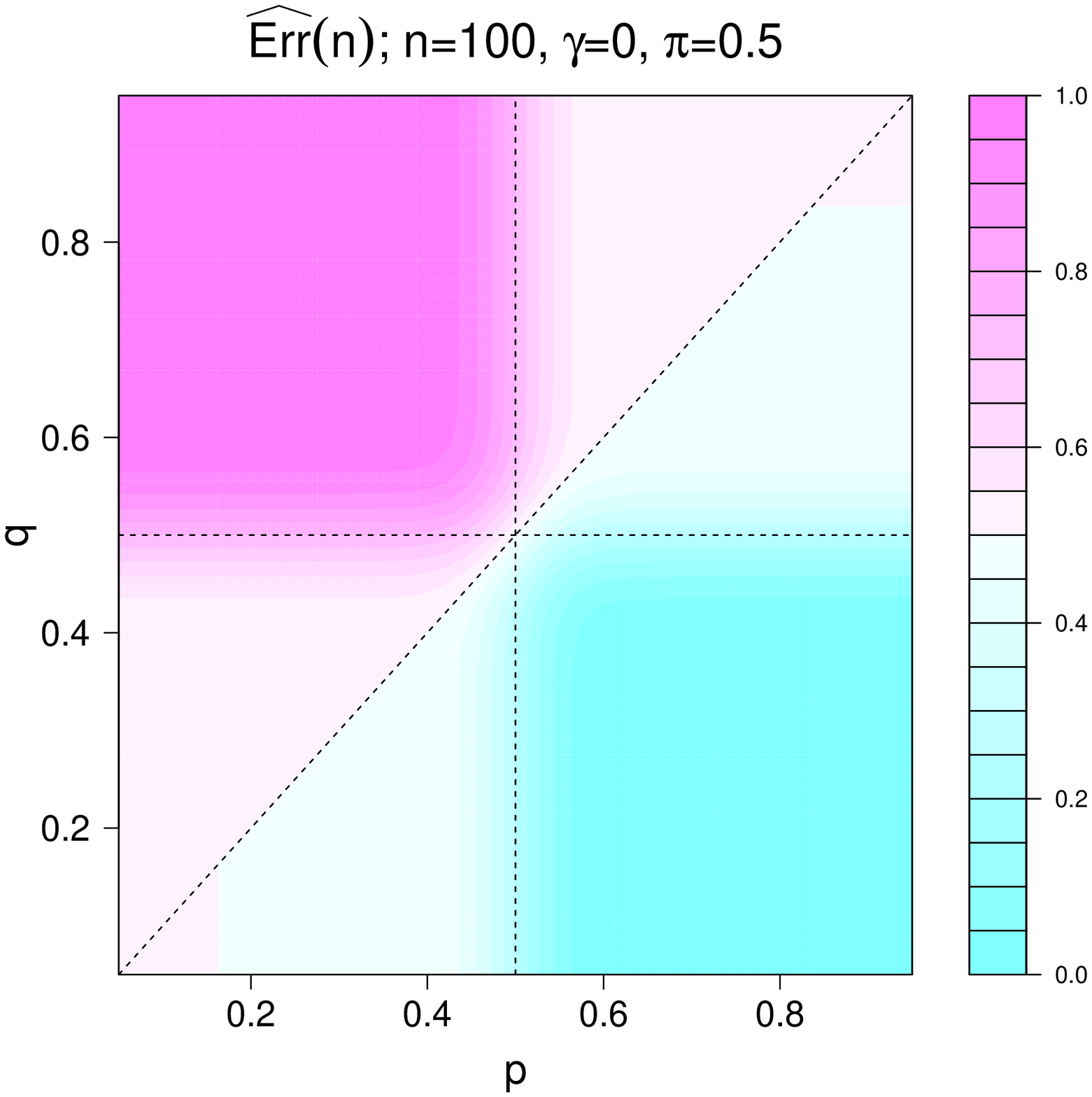}
\includegraphics[width=0.4\textwidth, angle=270]{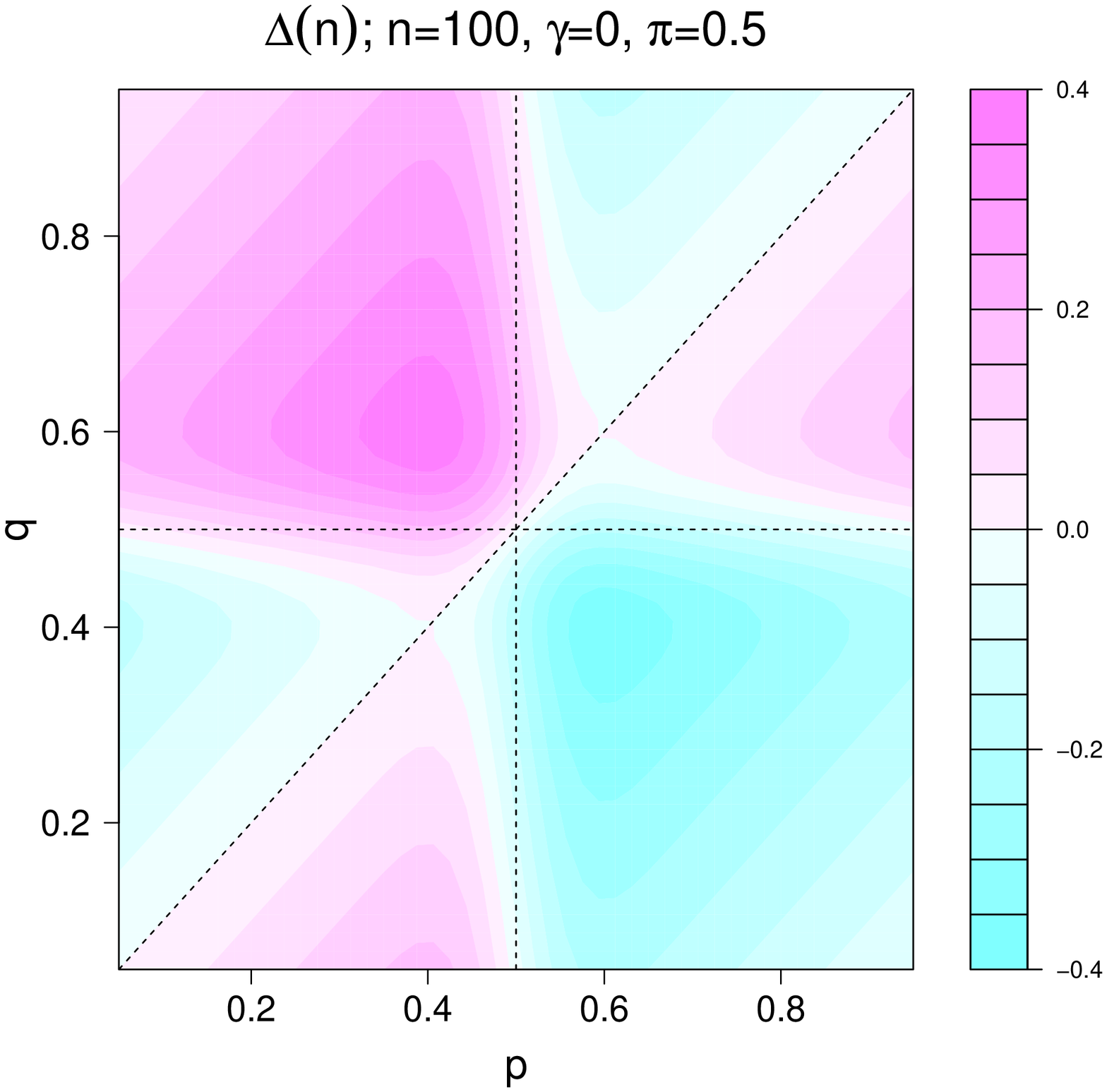} \\
\includegraphics[width=0.4\textwidth, angle=270]{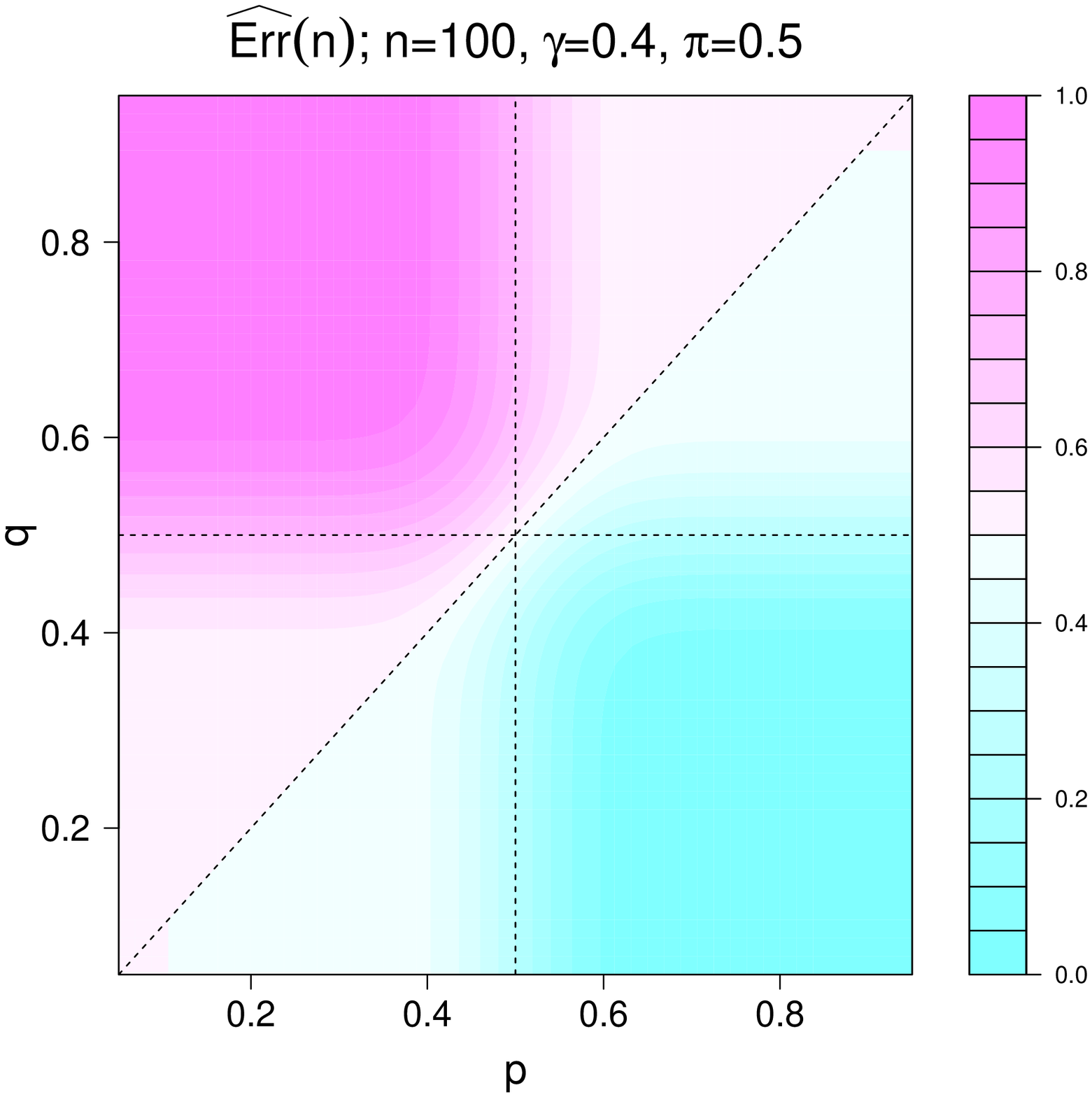}
\includegraphics[width=0.4\textwidth, angle=270]{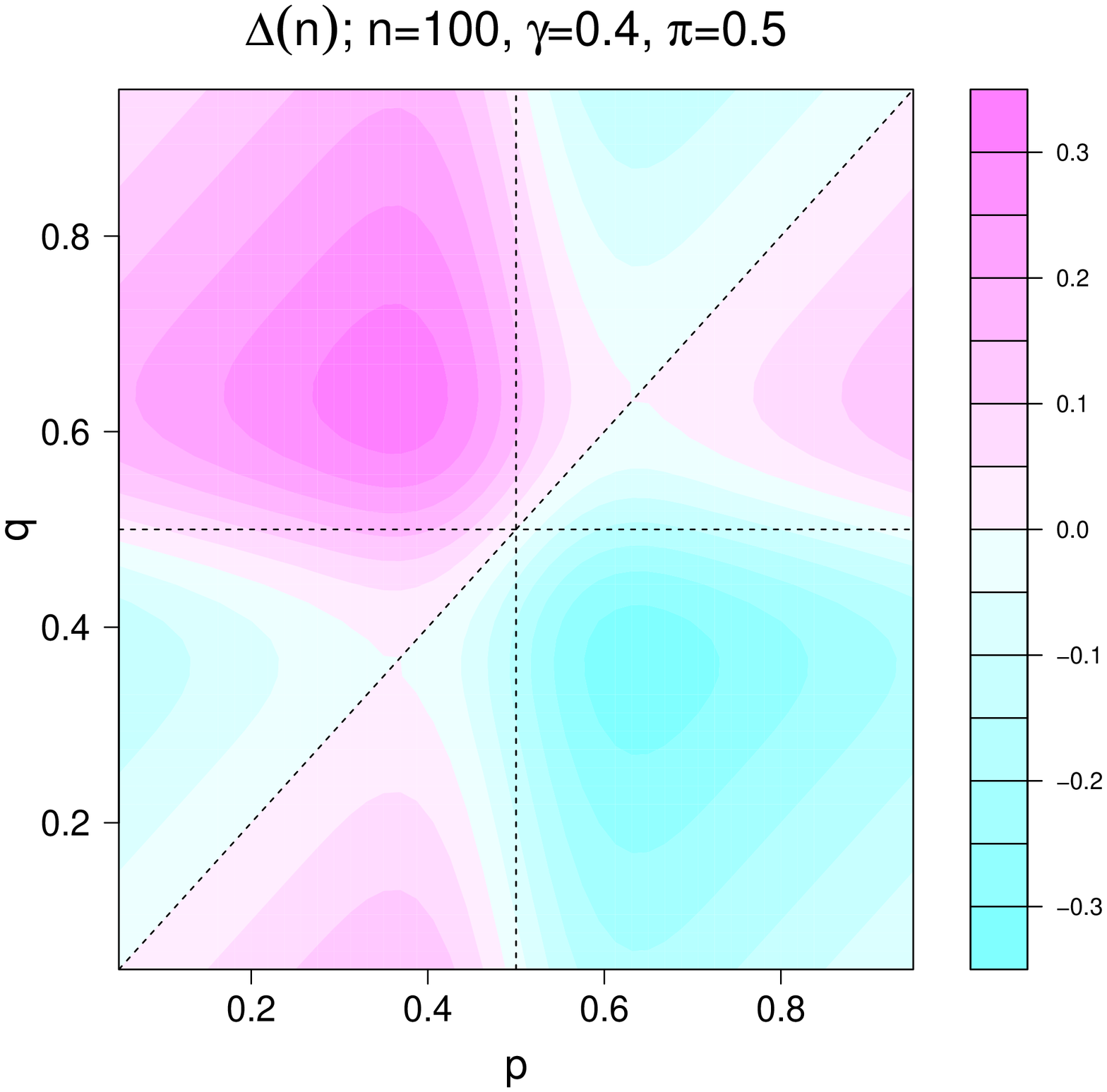} \\
\includegraphics[width=0.4\textwidth, angle=270]{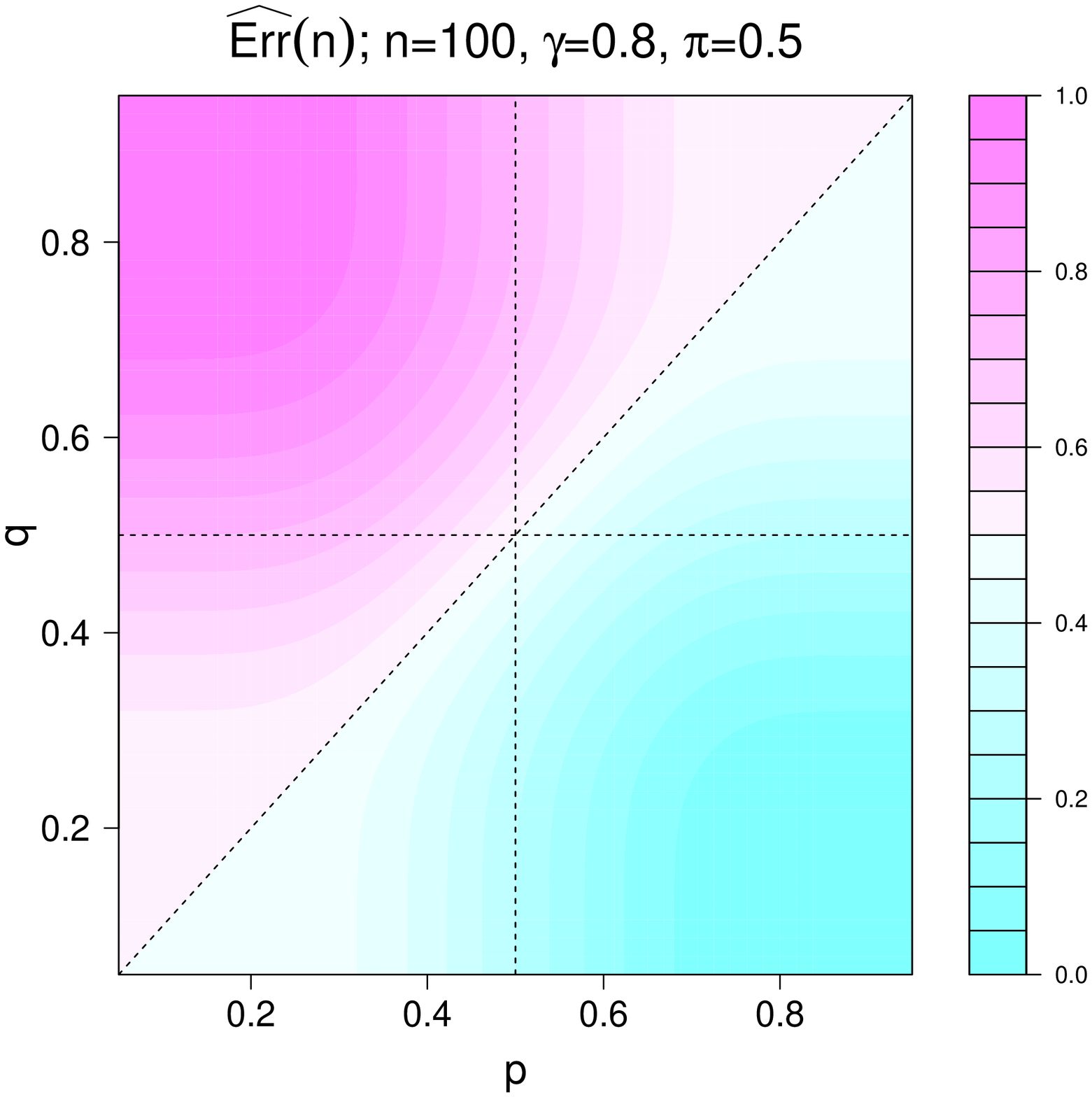}
\includegraphics[width=0.4\textwidth, angle=270]{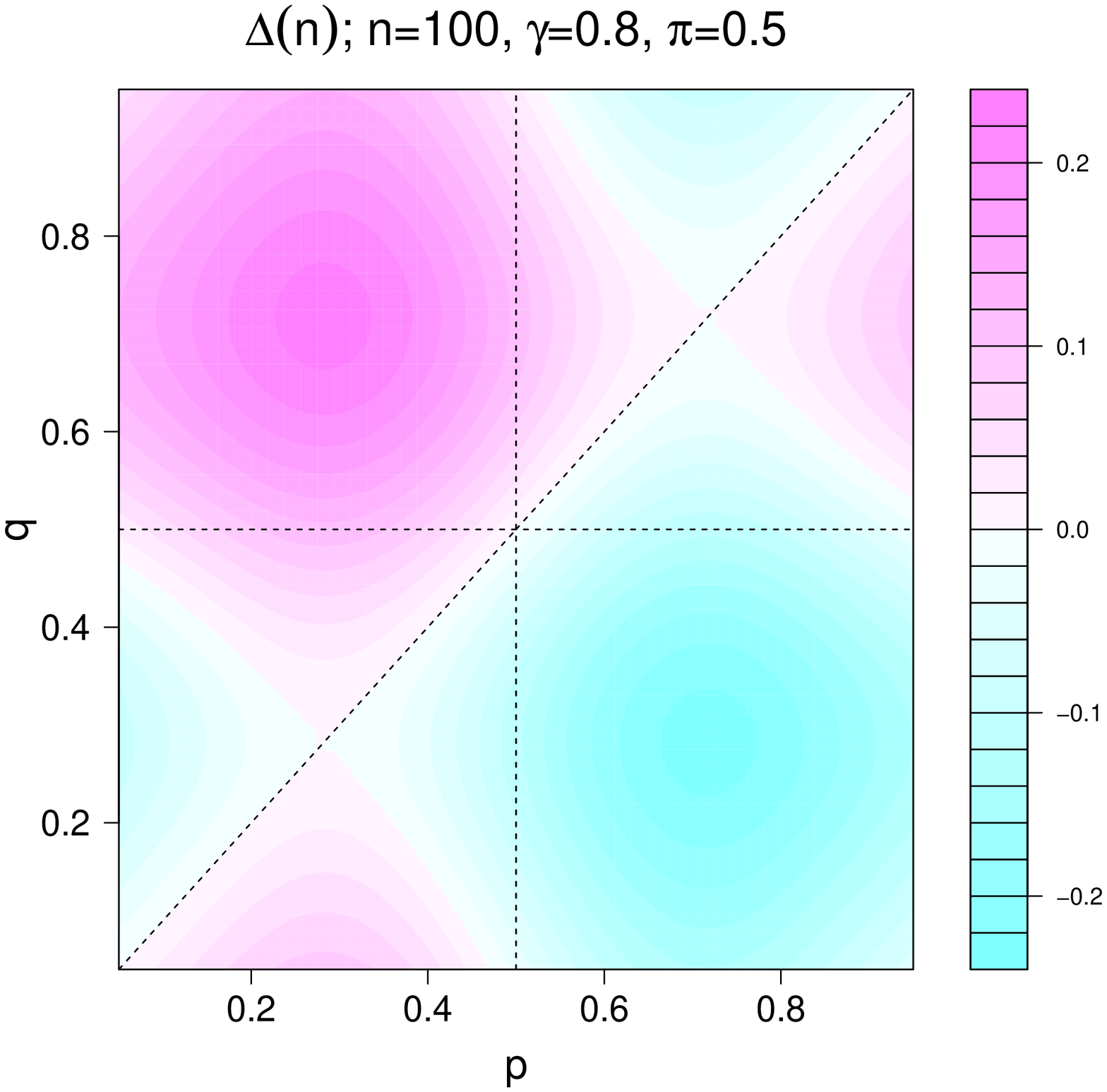} 
\caption{\label{fig:gamma}%
The case of $\corr_u(f_i,f_j)=\gamma^{|i-j|}$ for $u=p,q$ (Section~\ref{sec:AR}). 
This figure shows
$\widehat{\Err}(n)$ on the left and
$\Delta(n)$ on the right,
both as functions of $(p,q)$, for $\pi=1/2$, $n=100$,
and $\gamma=0$ (top), $\gamma=0.4$ (middle), $\gamma=0.8$ (bottom).}
\end{figure}

Figure~\ref{fig:gamma} shows the behaviors of $\widehat{\Err}(n)$ and 
$\Delta(n)$ for $\pi=1/2$, $n=100$, and $\gamma=0$ (uncorrelated case; 
Section~\ref{sec:indep}), $\gamma=0.4$, $\gamma=0.8$. We can see that the maximal 
amount of improvement obtainable by $g_n$ over a typical $f_i$ is about $40$ 
percentage points when $\gamma=0$, but only a little over $20$ percentage 
points when $\gamma=0.8$. In other words, when the classifiers are 
correlated, more of them are needed in order for $g_n$ to achieve the same level of 
performance --- hardly a surprising conclusion. Moreover, for the same 
$n=100$, both $\widehat{\Err}(n)$ and $\Delta(n)$ are clearly closer to 
their asymptotic limits (see Section~\ref{sec:main}) when $\gamma=0$. 
Overall, the effect of $\gamma>0$ is to slow down the convergence, 
although the same limits as in the (unrealistic) independent/uncorrelated 
case ($\gamma=0$) are eventually achieved.

Compared with the conditional independent case of Section~\ref{sec:indep}, 
the scenario considered in this section is slightly more 
realistic. The essence of the correlation model $\gamma^{|i-j|}$, of 
course, is that 
most of the correlations are kept quite small so the technical condition 
of Theorem~\ref{thm:main} can still be satisfied. 

\subsection{Case 3} 
\label{sec:equi}

For collections generated by an i.i.d.~mechanism such 
as (\ref{eq:RF}), it is perhaps even more realistic to consider the 
case that $\corr_u(f_i, f_j) = \lambda \in (0,1)$ for all $i \neq j$ and $u=p, q$. 
But then, we have
\beqn
\var_p(g_n) 
&=& \sum_{i=1}^n \var_p(f_i) + \sum_{i \neq j} \cov_p(f_i,f_j) \notag \\
&=& np(1-p) +  \lambda \times n(n-1)p(1-p) \notag\\
&=& n^2 \lambda p(1-p) + n(1-\lambda)p(1-p).  \label{eq:var-equi}
\eeqn
This means
\beqnn
\sigma_p^2 =
  \lim_{n\rightarrow\infty} \left[\frac{\var_p(g_n)}{n}\right]
= p(1-p) \lim_{n\rightarrow\infty} [ n \lambda + (1-\lambda) ]   
= \infty,
\eeqnn
and likewise for $\sigma_q^2$.
That is, 
the condition of Theorem~\ref{thm:main} no longer holds.

\subsection{The technical condition of Theorem~\ref{thm:main}}

Based on the three cases discussed above, it becomes fairly clear that the 
key technical condition of Theorem~\ref{thm:main} --- namely, $\sigma_p, \sigma_q < \infty$ --- is merely an indirect way to control the total amount of 
correlation among the individual classifiers in $\mathcal{F}_n$. In his ``random forest'' paper, \citet{randomForest} 
proved that, all other things being equal, classifiers having lower 
correlations with one another would form a better forest. Our analysis 
suggests that, even under fairly ideal conditions (well-controlled amount of 
total correlation), the majority-vote mechanism still requires more than 
$p>q$ in order to be truly beneficial.

\subsection{Case 3 (continued)}
\label{sec:equi-spec}

To get a sense of what could happen if the total amount of correlation 
should go ``out of control'', let us go back to Case 3 and ask a bold 
question: what if we went ahead and \emph{abused} the central limit theorem in 
this case, even though the usual requirement (\ref{eq:proper-condition}) is not met?
Since $\sigma_p, \sigma_q = \infty$, we replace $(\sigma_p\sqrt{n}, \sigma_q\sqrt{n})$ 
with $([\var_p(g_n)]^{1/2}, [\var_q(g_n)]^{1/2})$ in (\ref{eq:hatErr}). For large $n$,
(\ref{eq:var-equi}) implies that
\beqn
[\var_p(g_n)]^{1/2} \approx n \sqrt{\lambda p(1-p)} \quad\mbox{and}\quad
[\var_q(g_n)]^{1/2} \approx n \sqrt{\lambda q(1-q)}, \label{eq:var-equi-approx}
\eeqn
so our (abusive) estimate of $\Err(n)$ as $n\rightarrow\infty$ becomes
\beqn
\label{eq:hatErr-impropercase}
\widehat{\Err}(\infty) = 
  \left[ 
    \Phi\left( \frac{1/2-p}{\sqrt{\lambda p(1-p)}} \right)
  \right] \pi +
  \left[
  1-\Phi\left( \frac{1/2-q}{\sqrt{\lambda q(1-q)}} \right)
  \right] (1-\pi).
\eeqn
This is equivalent to the independent case 
(Section~\ref{sec:indep}) with $n=1/\lambda$. 
That is,
a large collection (of classifiers) 
with conditional pairwise correlation equal to $\lambda$ will behave
essentially like a conditionally independent (or uncorrelated) collection 
of size $1/\lambda$. Despite having abused the central limit theorem, this 
conclusion 
does not appear to be totally unreasonable after all.

\begin{figure}[htpb]
\centering
\includegraphics[width=0.4\textwidth, angle=270]{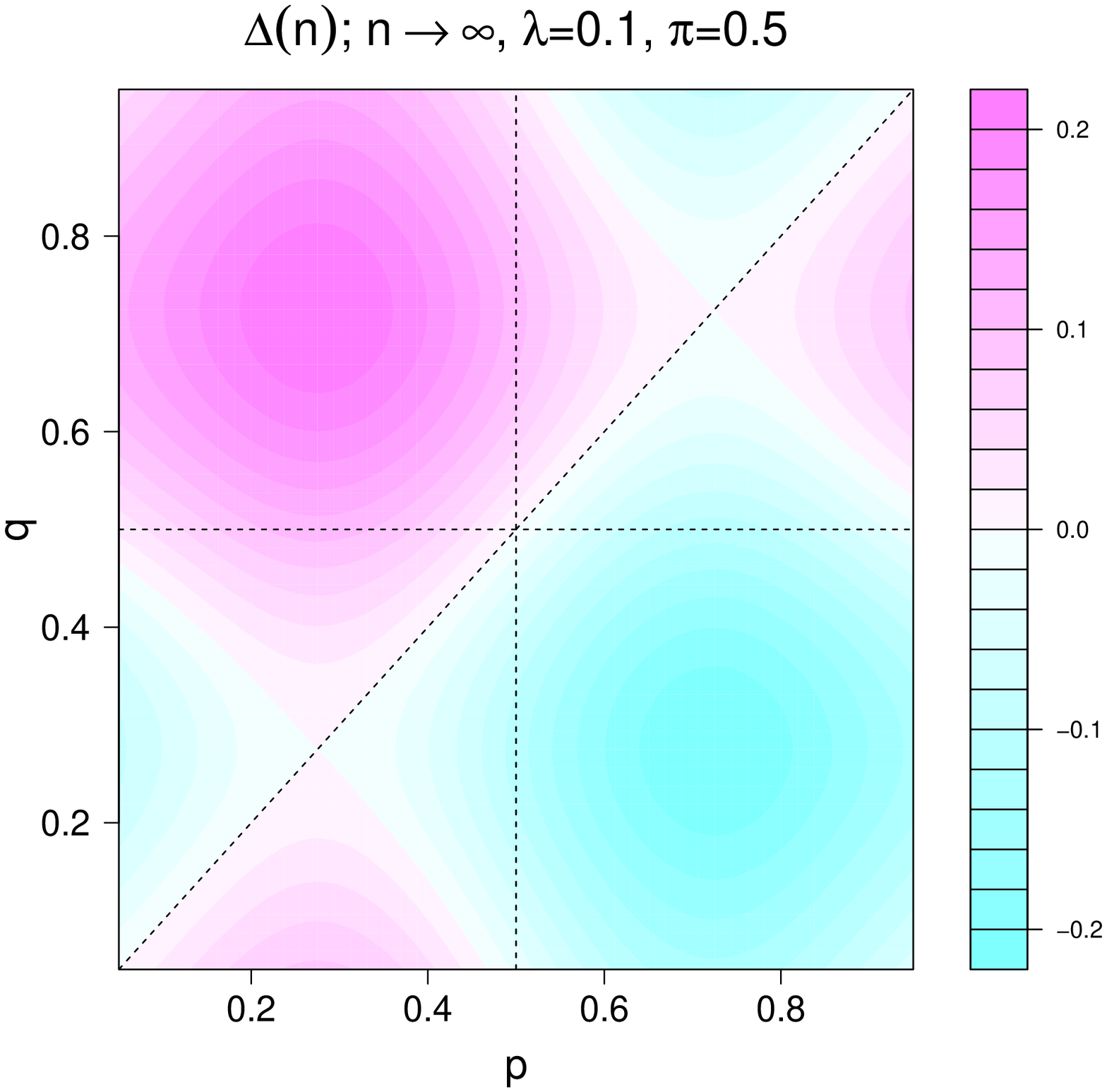}
\includegraphics[width=0.4\textwidth, angle=270]{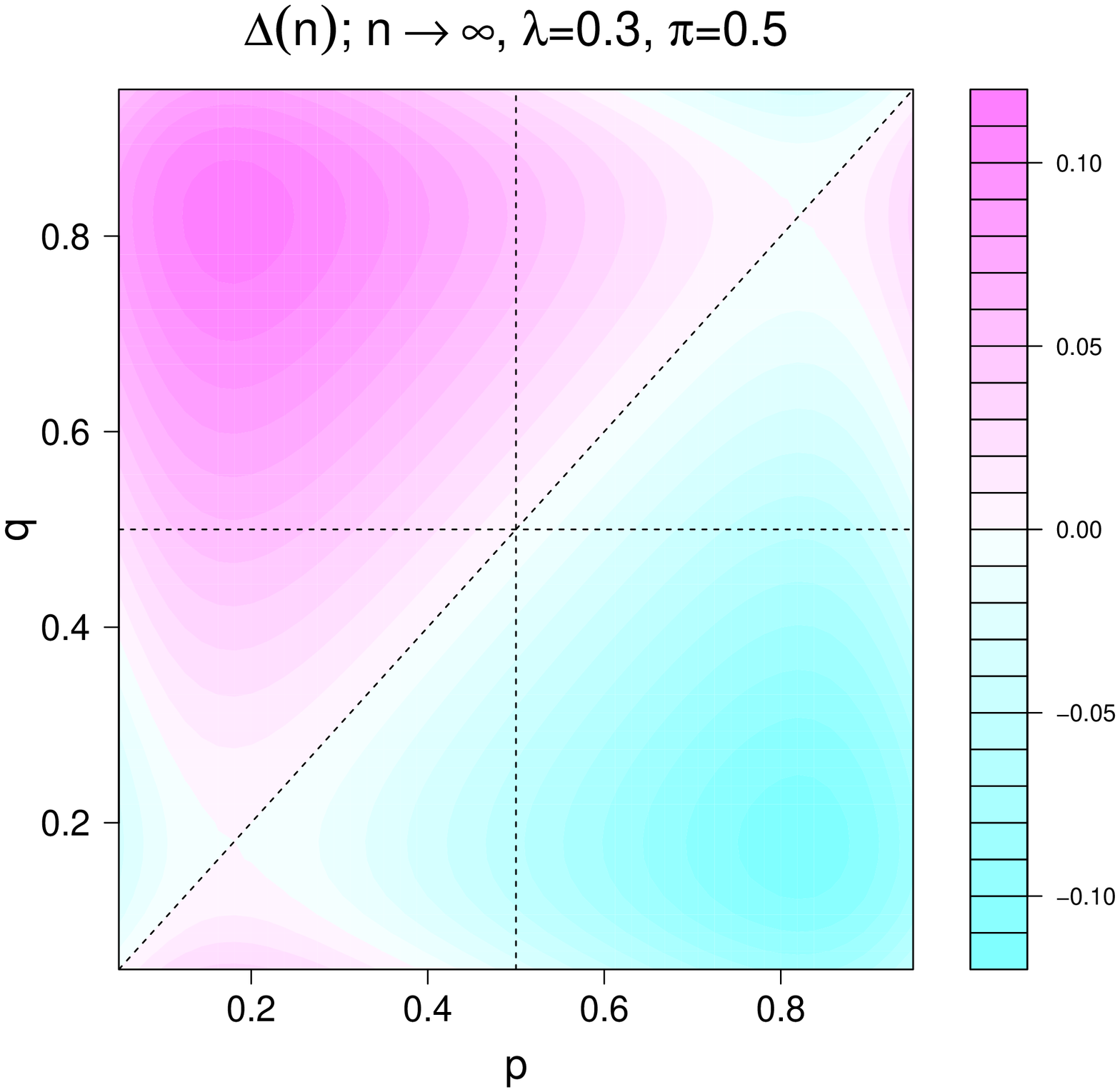} \\
\includegraphics[width=0.4\textwidth, angle=270]{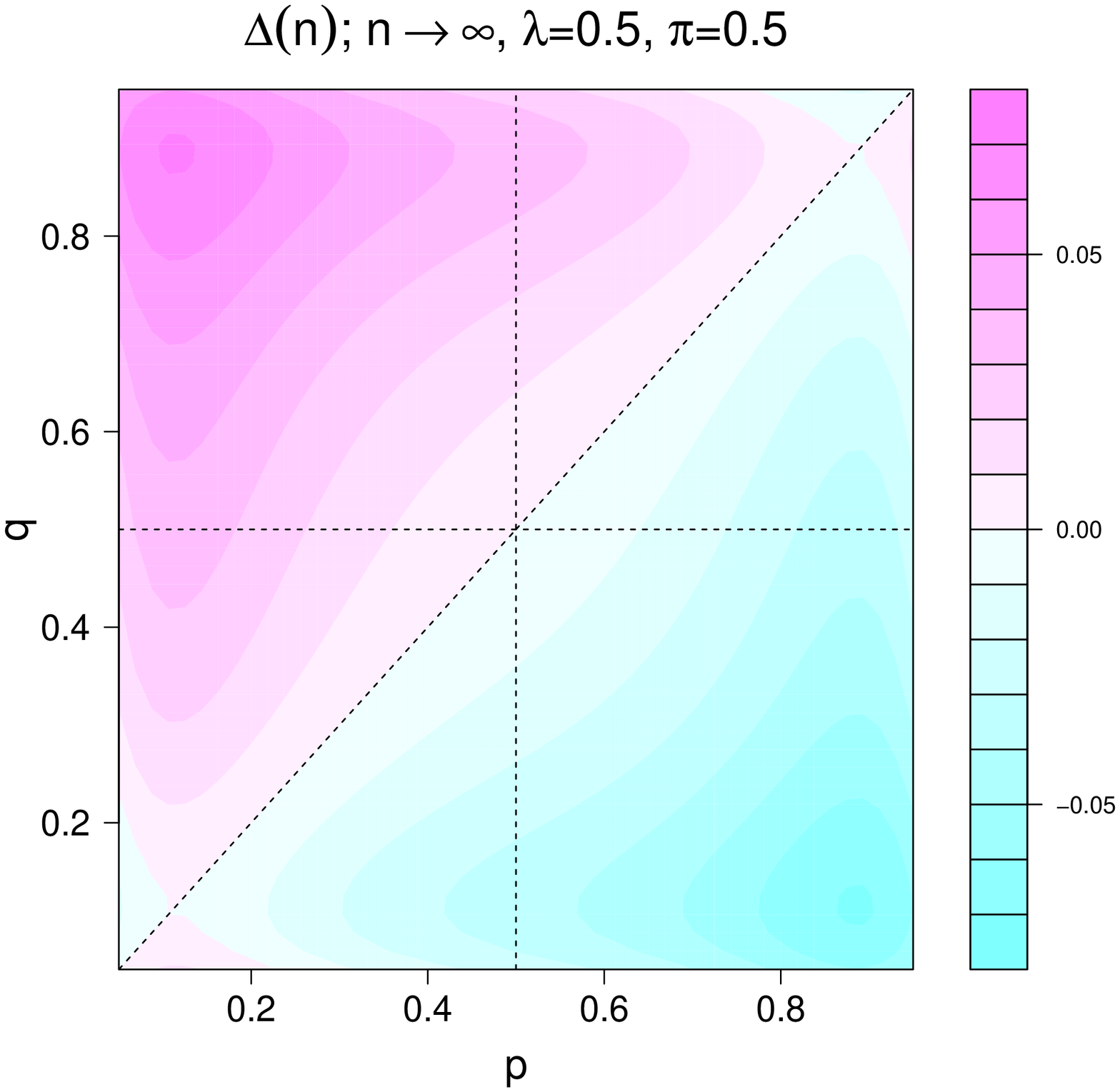}
\includegraphics[width=0.4\textwidth, angle=270]{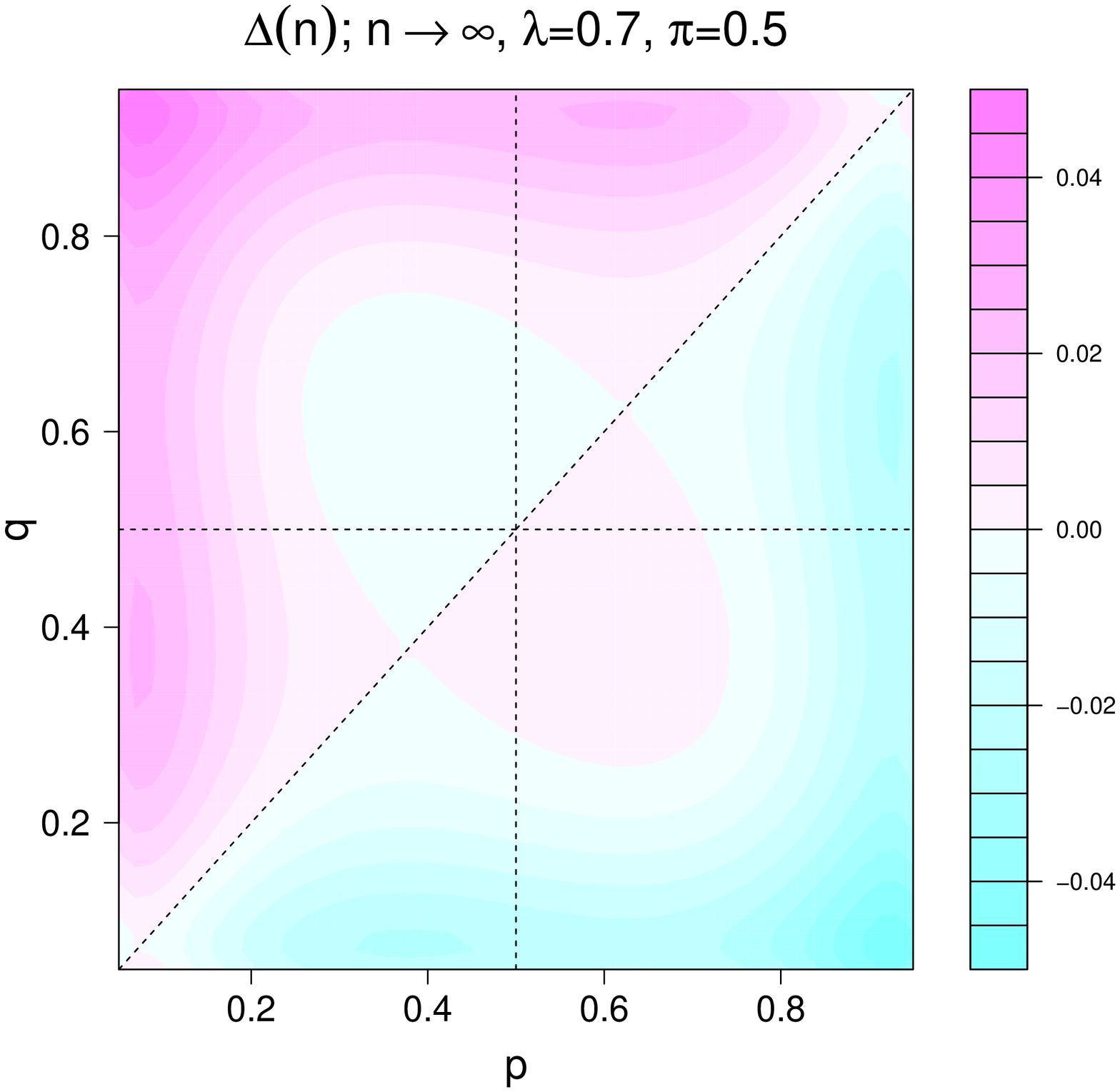}
\caption{\label{fig:lambda}%
The case of $\corr_u(f_i,f_j)=\lambda$ for $i \neq j$ and $u=p,q$ (Section~\ref{sec:equi-spec}).
This figure shows $\Delta(n)$ as a function of $(p,q)$, 
for $\pi=1/2$, $n\rightarrow\infty$, and $\lambda = 0.1, 0.3, 0.5, 0.7$. 
Notice that, since $\sigma_p,\sigma_q = \infty$ in this case, we have replaced $(\sigma_p\sqrt{n}, \sigma_q\sqrt{n})$ 
with $([\var_p(g_n)]^{1/2}, [\var_q(g_n)]^{1/2})$ in the expression of $\widehat{\Err}(n)$, and that 
any conclusions based on these plots are necessarily {\em speculative} as we can no longer assert $\Err(\infty)=\widehat{\Err}(\infty)$. }
\end{figure}

Figure~\ref{fig:lambda} shows the behavior of $\Delta(\infty)$ for $\pi=1/2$ 
and $\lambda=0.1, 0.3, 0.5, 0.7$, having replaced $(\sigma_p\sqrt{n}, \sigma_q\sqrt{n})$ in (\ref{eq:hatErr})
with $([\var_p(g_n)]^{1/2}, [\var_q(g_n)]^{1/2})$ as given in (\ref{eq:var-equi-approx}). 
Since the case of $\lambda=0.1$ 
(relatively low correlation) is 
``like'' having an independent collection of size $1/\lambda = 10$, the 
behavior of $\Delta(\infty)$ is still close enough to what our theory has 
predicted. For the case of $\lambda=0.7$ (very high correlation), however, 
even the ``golden 
region'' predicted by our theory --- namely $p \geq 1/2 \geq q$ --- fails 
to guarantee that $g_n$ is useful. The majority-vote mechanism {\em cannot} 
reduce the error rate {\em unless} the average FPR is very low or the average TPR is very high --- the south-east region shaped like a 
mirror image of the letter ``L''. Even then, the maximal amount of 
improvement obtainable is only a little over $4$ percentage points.

\section{Summary}

Ensemble classifiers can be constructed in many ways. \citet{weak-learner} 
showed that weak classifiers could be improved to achieve arbitrarily high 
accuracy, but never implied that a simple majority-vote mechanism could 
always do the trick. We have described an interesting phase transition phenomenon, which shows that, for the majority-vote 
mechanism to work, the weak classifiers cannot be ``too weak'' \emph{on average}. For example, in the case of equal priors, the collection must meet the minimum 
requirement of having an average TPR of at least $50\%$ and an average FPR of at most $50\%$, \emph{even when the 
total amount of correlation among the classifiers is well 
under control}. If 
the correlations are very high, this minimum requirement will likely have 
to be raised further, e.g., the classifiers may need to have either very 
high TPRs, or very low FPRs, or both 
(Section~\ref{sec:equi-spec}).